\documentclass[dvips,twoside]{article}
\usepackage[latin1]{inputenc}
\usepackage[T1]{fontenc}
\usepackage{parskip}
\usepackage{amsmath,amsfonts,amssymb,amsxtra}
\usepackage{latexsym}
\usepackage{theorem}
\usepackage{graphicx,psfrag}
\usepackage{verbatim}

\newtheorem{theo}{Theorem}
\newtheorem{prop}{Proposition}
\newtheorem{lemma}{Lemma}
\newtheorem{cor}{Corollary}

{\theoremstyle{plain}                       % Remark
\theorembodyfont{\rmfamily}
\newtheorem{Remark}{Remark}}

{\theoremstyle{plain}                       % Example
\theorembodyfont{\rmfamily}
\newtheorem{Example}{Example}}

{\theoremstyle{plain}                       % Definition
\theorembodyfont{\rmfamily}
\newtheorem{Definition}{Definition}}

\begin{document}

\pagestyle{myheadings} \markboth{\sc  I. Werner}{\sc The generalized
Markov measure as an equilibrium state}

\title{The generalized Markov measure as an equilibrium state.}
\author{Ivan Werner\footnote{The work was partially supported by EPSRC and School of
Mathematics and Statistics of University of St Andrews.}\\
   {\small Email: ivan\_werner@pochta.ru}}
\maketitle

\begin{abstract}\noindent
  In this paper, we continue development of the theory of contractive Markov
  systems (CMS) initiated in \cite{Wer1}. Also, this work can be seen as a small contribution to
  the theory of equilibrium states.

  We construct an energy function on the
  code space, using the coding map from \cite{Wer3}, and show that the generalized Markov measure associated with an irreducible
  CMS is a unique equilibrium state for this energy function
  if the vertex sets form an open partition of the state space of the CMS and the restrictions of the probability
  functions on their vertex sets are Dini-continuous and bounded away from
  zero.\\


\end{abstract}

\section{Introduction}%\renewcommand\baselinestretch 2
In  \cite{Wer1}, we introduced a theory of {\it contractive Markov
systems (CMS)} which provides a unifying framework in so-called
'fractal' geometry. It extends the known theory of {\it iterated
function systems (IFS) with place dependent probabilities}, which
are contractive on average, \cite{BDEG}\cite{Elton} in a way that it
also covers {\it graph directed constructions} of 'fractal' sets
\cite{MW}. In particular, Markov chains associated with such systems
naturally extend finite Markov chains and inherit some of their
properties.

By a {\it Markov system} we mean the following structure on a metric
space $(K,d)$, which generates a Markov process. It is given by a
family
\[\left(K_{i(e)},w_e,p_e\right)_{e\in E}\]
(see Fig. 1) where $E$ is the set of edges of a finite directed
(multi)graph $(V,E,i,t)$ ($V:=\{1,...,N\}$ is the set of vertices of
the directed (multi)graph (we do not exclude the case $N=1$),
$i:E\longrightarrow V$ is a map indicating the initial vertex of
each edge and $t:E\longrightarrow V$ is a map indicating the
terminal vertex of each edge), $K_1,K_2,...,K_N$ is a partition of
the metric space $K$ into non-empty Borel subsets, $(w_e)_{e\in E}$
is a family of Borel measurable self-maps on the metric space such
that $w_e\left(K_{i(e)}\right)\subset K_{t(e)}$ for all $e\in E$ and
$(p_e)_{e\in E}$ is a family of Borel measurable probability
functions on $K$ (i.e. $p_e(x)\geq 0$ for all $e\in E$ and
$\sum_{e\in E}p_e(x)=1$ for all $x\in K$) (associated with the maps)
such that each $p_e$ is zero on the complement of $K_{i(e)}$.

\begin{center}
\unitlength 1mm
\begin{picture}(70,70)\thicklines
\put(35,50){\circle{20}} \put(10,20){\framebox(15,15)}
\put(40,20){\line(2,3){10}} \put(40,20){\line(4,0){20}}
\put(50,35){\line(2,-3){10}} \put(5,15){$K_1$} \put(34,60){$K_2$}
\put(61,15){$K_3$} \put(31,50){\framebox(7.5,5)}
\put(33,45){\framebox(6.25,9.37)} \put(50,28){\circle{7.5}}
\put(45,21){\framebox(6,5)} \put(10,32.5){\line(6,1){15}}
\put(10,32.5){\line(3,-5){7.5}} \put(17.5,20){\line(1,2){7.5}}
\put(52,20){\line(2,3){4}} \put(13,44){$w_{e_1}$}
\put(35,38){$w_{e_2}$} \put(49,42){$w_{e_3}$}
\put(33,30.5){$w_{e_4}$} \put(30,15){$w_{e_5}$}
\put(65,37){$w_{e_6}$} \put(0,5){ Fig. 1. A Markov system}
\put(0,60){$N=3$} \thinlines \linethickness{0.1mm}
\bezier{300}(17,37)(20,46)(32,52)
\bezier{50}(32,52)(30.5,51.7)(30,49.5)
\bezier{50}(32,52)(30,51)(28.7,51.7)
\bezier{300}(26,31)(35,36)(35,47)
\bezier{50}(35,47)(35,44.5)(33.5,44)
\bezier{50}(35,47)(35,44)(36,44) \bezier{300}(43,50)(49,42)(51,30)
\bezier{50}(51,30)(50.5,32)(49.2,32.6)
\bezier{50}(51,30)(50.6,32)(51.5,33.2)
\bezier{300}(39,20)(26,17)(18,25)
\bezier{50}(18,25)(19.5,24)(20,21.55)
\bezier{50}(18,25)(20,23.5)(22,24) \bezier{300}(26,26)(37,28)(47,24)
\bezier{50}(47,24)(45,25)(43,24) \bezier{50}(47,24)(45,25)(44,26.5)
\bezier{100}(54.5,31.9)(56,37.3)(61,36.9)
\bezier{100}(61,36.9)(64.5,36.5)(66,34)
\bezier{100}(66,34)(68,30.5)(64.9,26.8)
\bezier{100}(64.9,26.8)(61.6,23.3)(57,23)
\bezier{50}(57,23)(58.5,23.3)(60.1,22.7)
\bezier{50}(57,23)(58.8,23.3)(59.5,24.8)
\end{picture}
\end{center}

A Markov system is called {\it irreducible} or {\it aperiodic} iff
its directed graph is irreducible or aperiodic respectively. We call
a Markov system $\left(K_{i(e)},w_e,p_e\right)_{e\in E}$ {\it
contractive} with an {\it average contracting rate} $0<a<1$ iff it
satisfies the following {\it condition of contractiveness on
average}:
\begin{equation}\label{cc}
 \sum\limits_{e\in E}p_e(x)d(w_ex,w_ey)\leq ad(x,y)\mbox{ for all
}x,y\in K_i,\ i=1,...,N.
\end{equation}
This condition was discovered by R. Isaac in 1961 \cite{Is} (for the
case $N=1$).

Markov system $\left(K_{i(e)},w_e,p_e\right)_{e\in E}$ determines a
Markov operator $U$ on the set of all bounded Borel measurable
functions $\mathcal{L}^0(K)$ by
\[Uf:=\sum\limits_{e\in E}p_ef\circ w_e\mbox{ for all
 }f\in\mathcal{L}^0(K)\] and its adjoint operator $U^*$ on the set of all Borel probability
 measures $P(K)$ by
\[U^*\nu(f):=\int U(f)d\nu\mbox{ for all }f\in\mathcal{L}^0(K)\mbox{ and }\nu\in P(K).\]

\begin{Remark}
    Note that each map $w_e$ and each probability $p_e$ need to
    be defined only on the corresponding vertex set $K_{i(e)}$. This
   is sufficient for the condition (\ref{cc}) and the definition of
   $U^*$. For the definition of $U$, we can consider each $w_e$ to
   be extended on the whole space $K$ arbitrarily and each $p_e$ to
   be extended on $K$ by zero.

   Also, the situation applies where each vertex set $K_i$ has its
   own metric $d_i$. In this case, one can set
   \[d(x,y)=\left\{\begin{array}{cc}
       d_i(x,y) & \mbox{ if }x,y\in K_i \\
       \infty & \mbox{otherwise}
     \end{array}\right.
   \] and use the convention $0\times\infty=0$.
\end{Remark}

 We say $\mu\in P(K)$ is an invariant measure of the CMS iff
$U^*\mu=\mu$. A Borel probability measure $\mu$ is called {\it
 attractive} measure of the CMS if
 \[{U^*}^n\nu\stackrel{w^*}{\to}\mu\mbox{ for all }\nu\in P(K),\]
 where $w^*$ means weak$^*$ convergence.
Note that an attractive probability measure is a unique invariant
probability measure of the CMS if $U$ maps continuous functions on
continuous functions. The following theorem states some properties
of such CMSs.

\begin{theo}\label{ft} Let $K$ be a metric space in which sets of
finite diameter are relatively compact.  Suppose
$\left(K_{i(e)},w_e,p_e\right)_{e\in E}$ is a CMS with an average
contracting rate $0<a<1$ such that the family $K_1,...,K_N$
partitions $K$ into open subsets and each $p_e|_{K_{i(e)}}$ is
continuous on $K_{i(e)}$. Then the following hold:\\
(i) The sequence $\left({U^*}^k\delta_x\right)_{k\in\mathbb{N}}$ is
tight for all $x\in K$, i.e. for all $\epsilon>0$, there exists a
compact subset $Q\subset K$
such that ${U^*}^k\delta_x(Q)\geq 1-\epsilon$ for all $k\in\mathbb{N}$.\\
(ii) The CMS has an invariant Borel probability measure $\mu$.\\
(iii) The invariant probability measure $\mu$ is unique iff
 \[\frac{1}{n}\sum\limits_{k=1}^nU^kg(x)\to\int gd\mu\mbox{ for
 all }x\in K\mbox{ and }g\in C_B(K).\]
(iv) If the invariant probability measure is unique, then
\[\sum\limits_{i=1}^N\int\limits_{K_i}d(x,x_i)d\mu(x)<\infty\mbox{
 for all }x_i\in K_i,\ i=1,...,N.\]
\end{theo}
For the proof see \cite{Wer1}.

Furthermore, it was shown in \cite{Wer1} that contractive Markov
systems inherit some properties of finite Markov chains if the
restrictions of the probabilities on their vertex sets satisfy some
stronger conditions. Namely, it has a unique invariant probability
measure in an irreducible case and an attractive probability measure
in an aperiodic case if the vertex sets $K_1,...,K_N$ form an open
partition of the state space and each $p_e|_{K_{i(e)}}$ is
Dini-continuous and bounded away from zero (see also \cite{Wer2} for
the case of constant probabilities $p_e|_{K_{i(e)}}$ and compact
state space).

A function $h:(X,d)\longrightarrow\mathbb{R}$ is called {\it
Dini-continuous} iff
 for some $c >0$  \[\int_0^c\frac{\phi(t)}{t}dt<\infty\]
  where $\phi$ is {\it the modulus of uniform continuity} of $h$, i.e.
     \[\phi(t):=\sup\{|h(x)-h(y)|:d(x,y)\leq t,\ x,y\in X\}.\]
It is easily seen that Dini-continuity is weaker than H\"{o}lder and
stronger than the uniform continuity.  There  is a well known
characterization of Dini-continuity, which will be useful later.
\begin{lemma}\label{Dc}
 Let $0<c<1$ and $b>0$.
  A  function $h$ is Dini-continuous  iff
 \[\sum_{n=0}^\infty\phi\left(bc^n\right)<\infty\] where $\phi$ is
 the modulus of uniform continuity of $h$.
\end{lemma}
The proof is simple (e.g. see \cite{Wer1}).

Also, associated with the Markov system is a measure preserving
transformation $S:(\Sigma,\mathcal{B}(\Sigma),
M)\longrightarrow(\Sigma,\mathcal{B}(\Sigma), M)$, which we call a
{\it generalized Markov shift}, where
$\Sigma:=\{(...,\sigma_{-1},\sigma_0,\sigma_1,...):\sigma_i\in E\
\forall i\in\mathbb{Z}\}$ is the {\it code space} provided with the
product topology, $\mathcal{B}(\Sigma)$ denotes Borel
$\sigma$-algebra on $\Sigma$ and $M$ is a {\it generalized Markov
measure} on $\mathcal{B}(\Sigma)$ given by
\begin{equation*}\label{gMm}
   M\left(_m[e_1,...,e_k]\right):=\int
p_{e_1}(x)p_{e_2}(w_{e_1}x)...p_{e_k}(w_{e_{k-1}}\circ...\circ
w_{e_1}x)d\mu(x)
\end{equation*} for every cylinder set
$_m[e_1,...,e_k]:=\{\sigma\in\Sigma:\
\sigma_m=e_1,...,\sigma_{m+k-1}=e_k\}$, $m\in\mathbb{Z}$, where
$\mu$ is an invariant Borel probability measure of the Markov
system, and $S$ is the usual left shift map on $\Sigma$. It is easy
to verify that $S$ preserves measure $M$, since $U^*\mu=\mu$ (see
\cite{Wer3}).

For a contractive Markov system (CMS), the Markov process associated
with the CMS can be represented as a factor of the generalized
Markov shift via a {\it coding map} $F:(\Sigma,\mathcal{B}(\Sigma),
M)\longrightarrow K$ which was constructed in \cite{Wer3}. It is
defined by
\[F_{x_1,...,x_N}(\sigma):=\lim\limits_{m\to-\infty}w_{\sigma_0}\circ
w_{\sigma_{-1}}\circ...\circ w_{\sigma_m}x_{i(\sigma_m)}\mbox{ for
}M\mbox{-a.e. }\sigma\in\Sigma,\] under some conditions (see the
next lemma), where $x_i\in K_i$ for each $i=1,...,N$.

\begin{lemma}\label{cm}
Suppose that  $p_e|_{K_{i(e)}}$ is Dini-continuous and there exists
$\delta>0$ such that $p_e|_{K_{i(e)}}\geq\delta$ for all $e\in E$.
 Let $x_i,y_i\in K_i$ for all  $1\leq i\leq N$. Then the following hold:

 (i)  $F_{x_1...x_N}$ is defined $M$-a.e.,

 (ii) $F_{x_1...x_N}=F_{y_1...y_N}$ $M$-a.e.. and

 (iii) There exists a sequence of closed subsets $Q_1\subset
 Q_2\subset...\subset\Sigma$ with\\ $\lim_{k\to\infty}M(Q_k)=1$
 such that all $F_{x_1...x_N}|_{Q_k}$ are locally
 H\"{o}lder-continuous with the same H\"{o}lder-constants.
\end{lemma}
In the following, we fix $x_1,...,x_N$ and denote the coding map
simply by $F$. This coding map is the key tool for our
investigation.

\begin{Example}[decimal expansion]
 Consider ten maps $w_e$, $e\in E:=\{0,...,9\}$, on $([0,1], |.|)$ given by
 $w_e(x):=1/10x+e/10$ for all $x\in [0,1]$. Obviously, for any
 family of probability functions $p_e$, $e\in E$, the family
 $([0,1],w_e,p_e)_{e\in E}$ is a CMS. The coding map for this CMS is nothing else as the
usual decimal expansion of real numbers from $[0,1]$.
\end{Example}

\begin{Example}
  Let $G:=(V,E,i,t)$ be a finite irreducible directed
  (multi)graph. Let $\Sigma^-_G:=\{(...,\sigma_{-1},\sigma_0):\
  \sigma_m\in E\mbox{ and }
 t(\sigma_m)=i(\sigma_{m-1})\ \forall
 m\in\mathbb{Z}\setminus\mathbb{N}\}$ ({\it one-sided
  subshift of finite type} associated with $G$) endowed with the
  metric $d(\sigma,\sigma'):=2^k$ where $k$ is the smallest
  integer with $\sigma_i=\sigma'_i$ for all $k<i\leq 0$. Let $g$
  be a positive, Dini-continuous function on $\Sigma_G$ such that
  \[\sum\limits_{y\in T^{-1}(\{x\})}g(y)=1\mbox{ for all }x\in\Sigma_G\]
  where $T$ is the right shift map on $\Sigma^-_G$. Set
  $K_i:=\left\{\sigma\in\Sigma^-_G:t(\sigma_0)=i\right\}$  for
  every $i\in V$ and, for
  every $e\in E$,
  \[w_e(\sigma):=(...,\sigma_{-1},\sigma_{0},e),\ p_e(\sigma):=g(...,\sigma_{-1},\sigma_{0},e)
  \mbox{ for all }\sigma\in K_{i(e)}.\]
  Obviously, maps $(w_e)_{e\in E}$ are contractions. Therefore,
  $\left(K_{i(e)}, w_e, p_e\right)_{e\in E}$
   defines a CMS. An invariant probability measure of such a CMS is called a $g$-measure. This
  notion was introduced by M. Keane \cite{Ke}.
  In this example, the coding map $F:\Sigma_G\longrightarrow\ \Sigma^-_G$
  is nothing else as the natural projection, and the generalized Markov measure is the natural
  extension of the $g$-measure (or in other words, the $g$-measure is the image of the generalized
  Markov measure under the coding map). The equivalence between $g$-measures and
  equilibrium states for $\log g$ was elaborated
  by Ledrappier \cite{Le}. See also \cite{Wal} for a proof of the uniqueness of the
  equilibrium  state in this example.
\end{Example}

In this paper, we are concerned with the following question. Is the
generalized Markov measure associated with a CMS a unique
equilibrium state for some energy function?  We report here that the
answer to this question is YES, under some conditions (see Corollary
\ref{DCc}). The result seems to be beyond the well known
thermodynamic formalism. It turns out that our energy function is
not upper semicontinuous in general (see Example \ref{ex}). Recall
that the upper semicontinuity of an energy function is a widely used
condition in the rigorous theory of equilibrium states which insures
that the convex set of equilibrium states is non-empty and compact
in the weak$^*$ topology (see e.g. \cite{Kel}).

Also, an interesting point of the presented result is that it
introduces some symbolic dynamical systems of infinite type, which
have the full measures $M$.

\begin{Definition}
Let $X$ be a metric space and $T$ a continuous transformation on it.
Denote by $P(X)$ the set of all
 Borel probability measures on $X$  and by $P_T(X)$
the set of all  $T$-invariant Borel probability measures on $X$.
 We call a Borel measurable
 function $f :X \longrightarrow [-\infty,0]$ an {\it energy function}. Suppose that $T$ has a finite
 topological entropy, i.e. $\sup_{\Theta\in
P_T(X)}h_\Theta(T)<\infty$, where $h_{\Theta}(T)$ is the
Kolmogorov-Sinai entropy of $T$ with respect to measure $\Theta$. We
call
 \[P(f)=\sup\limits_{\Theta\in
P_T(X)}\left(h_\Theta(T)+\Theta(f)\right)\]  the {\it pressure} of
$f$.  We call $\Lambda\in P_T(X)$ {\it an equilibrium state} for $f$
iff
\[h_{\Lambda}(T)+\Lambda(f)=P(f).\]
\end{Definition}
Note that $\sup_{\Theta\in P_S(\Sigma)}h_\Theta(S)=\log|E|$ (e.g.
Example 4.2.6 in \cite{Kel}).

\section{Main part}

Let
\begin{equation*}
    \mathcal{M}:=\left(K_{i(e)},w_e,p_e\right)_{e\in E}
\end{equation*} be a contractive Markov
system with an average contracting rate $0<a<1$ and an invariant
Borel probability measure $\mu$. We assume that: $(K,d)$ is a metric
space in which sets of finite diameter are relatively compact and
the family $K_1,...,K_N$ partitions $K$ into non-empty open subsets;
each probability function $p_e|_{K_{i(e)}}$ is uniformly continuous
and bounded away from zero by $\delta>0$; the set of edges $E$ is
finite and the map $i:E\longrightarrow V$ is surjective. Note that
the assumption on the metric space implies that it is locally
compact separable and complete. We shall denote the space of all
bounded continuous functions on $K$ by $C_B(K)$.

We consider $\Sigma$ endowed with the metric
$d'(\sigma,\sigma'):=(1/2)^k$ where $k$ is the largest integer with
$\sigma_i=\sigma'_i$ for all $|i|<k$. Fix $x_i\in K_{i}$ for all
$i=1,...,N$.

The construction of the energy function goes through a definition of
an appropriate shift invariant subset of $\Sigma$ on which the
energy function shall be finite.

 Let
\[\Sigma_G:=\{\sigma\in\Sigma:\ t(\sigma_j)=i(\sigma_{j+1})\ \forall
j\in\mathbb{Z}\},\]
\begin{equation*}
    D:=\{\sigma\in\Sigma_G:\
\lim\limits_{m\to-\infty}w_{\sigma_0}\circ
w_{\sigma_{-1}}\circ...\circ w_{\sigma_m}x_{i(\sigma_m)}\mbox{
exists}\}
\end{equation*} and
\[Y:=\bigcap\limits_{i=-\infty}^\infty S^i(D).\]

\begin{lemma}\label{pY}
(i) $F(\sigma)$ is defined for all $\sigma\in Y$ and $Y$ is a shift invariant subset of $\Sigma_G$.\\
(ii) If each $p_e|_{K_{i(e)}}$ is Dini-continuous, then $M(Y)=1$.\\
(iii) If each $p_e|_{K_{i(e)}}$ is Dini-continuous and the CMS has
an invariant probability measure $\mu$ such that $\mu(K_{i(e)})>0$
for all $e\in E$. Then $Y$ is dense in $\Sigma_G$.
\end{lemma}

{\it Proof.} (i) is clear, by the definitions of $F$ and $Y$.
 Note that the condition of the contractiveness on average and the boundedness away from zero of the
 functions $p_e|_{K_{i(e)}}$ imply that each map $w_e|_{K_{i(e)}}$ is continuous (Lipschitz).
 Therefore, $S(D)\subset D$. By Lemma \ref{cm}, $M(D)=1$. This implies that $M(Y)=1$. If in addition
$\mu(K_{i(e)})>0$ for all $e\in E$, then  $M(O)>0$ for every open
$O\subset\Sigma_G$. This implies  $(iii)$. \hfill$\Box$

\begin{Remark}
(i) Note that $Y$ and $F$ depend on the choice of $x_i$'s. By
Corollary 1 in \cite{Wer3},
$Y$ changes only modulo $M$-zero set by a different choice of $x_i$'s.  \\
(ii) If all maps $w_e|_{K_{i(e)}}$ are contractive, then
$Y=\Sigma_G$ and $F|_{\Sigma_G}$ is H\"{o}lder-continuous (easy to
check).
\end{Remark}

In the following, we are going to represent the generalized Markov
measure $M$ as a unique equilibrium state for the energy function
$u$ given by
\begin{equation}\label{ef}
    u(\sigma)=\left\{\begin{array}{cc}
    \log p_{\sigma_1}\circ F(\sigma)&  \mbox{if }\sigma\in Y \\
    -\infty& \mbox{if }\sigma\in\Sigma\setminus Y
  \end{array}\right.
\end{equation}
if the CMS has a unique invariant Borel probability measure and
$M(Y)=1$.

Let's consider a simple example which illustrates $F$, $Y$ and $u$.
\begin{Example}\label{ex}
 Let $(K,d)=(\mathbb{R},|.|)$. Consider two maps
 \[w_0(x):=\frac{1}{2}x,\ w_1(x):=2x\mbox{ for all }x\in\mathbb{R}\]
 with probability functions
 \[p_0(x):=\frac{1}{6}\sin^2x+\frac{17}{24},\ p_1(x):=\frac{1}{6}\cos^2x+\frac{1}{8}\mbox{
 for all }x\in\mathbb{R}.\]
  Then a simple
 calculation shows that $(\mathbb{R},w_e,p_e)_{e\in\{0,1\}}$
 defines a CMS with an average contracting rate $45/48$. In this
 case, $\Sigma_G=\{0,1\}^\mathbb{Z}$. If we take $x=0$ for the definition of $Y$,
 then, obviously, $Y=\Sigma_G$. Now, let $x\neq 0$. Let $N_{0n}(\sigma)$ and $N_{1n}(\sigma)$
 be the numbers of zeros and ones in $(\sigma_{-n},...,\sigma_0)$
 respectively for every $\sigma\in\Sigma_G$. Then, obviously,
 $\sigma\notin Y$ if
 $\left(N_{1n}(\sigma)-N_{0n}(\sigma)\right)\to\infty$. Hence,
 $Y\neq\Sigma_G$ and, by Lemma \ref{pY} (iii), $Y$ is a
 dense shift invariant subset of $\Sigma_G$. Since $Y$ is not closed, $u$ is not upper
 semicontinuous.

 Also, it is not difficult to see that in a general case there is
 no hope to find $x_i$ such that $u$ becomes upper
 semicontinuous, e.g. change $w_1$ to $w_1(x)=2x+1$, then, for any
 choice of $x$ for the definition of $Y$, $Y\neq\Sigma_G$.
\end{Example}

In what follows, we shall denote the restrictions of $F$ and $M$ to
$Y$ by the same letters.

\begin{Definition}Let $\mathcal{A}_0$ be the finite $\sigma$-algebra on $\Sigma$
generated by the zero-time partition. Set
$\mathcal{F}:=\bigvee_{i=0}^\infty S^i\mathcal{A}_0$.
 Let $\Lambda\in P_S(\Sigma)$. Define, for $e\in E$,
  \[g_{\Lambda e}:=E_{\Lambda}\left(1_{_1[e]}|\mathcal{F}\right)\mbox{ and } g_\Lambda(\sigma):=
  g_{\Lambda{\sigma_1}}(\sigma)\mbox{ for }\Lambda\mbox{-a.e. }\sigma\in \Sigma,\]
 where $E_\Lambda(.|.)$ denotes the conditional expectation with
 respect to $\Lambda$. Note that we can consider $g_{\Lambda
 e}(\sigma)=g_{\Lambda e}(...,\sigma_{-1},\sigma_0)$ for all $\sigma\in\Sigma$.
\end{Definition}

We are going to show that every $\Lambda\in P_S(\Sigma)$ is an
equilibrium state for a version of $\log g_\Lambda$, but first let's
see some properties of the function $g_\Lambda$.

\begin{lemma}\label{ef}
Let $\Lambda\in P_S(\Sigma)$.\\
 $(i)$ $0\leq g_\Lambda\leq 1$ $\Lambda$-a.e..\\
 $(ii)$ $\sum_{e\in E}g_{\Lambda e}=1$ $\Lambda$-a.e..\\
 $(iii)$ $\Lambda\left(\{g_\Lambda=0\}\right)=0$.\\
 $(iv)$ $1/g_\Lambda,\ \log g_\Lambda\in \mathcal{L}^1(\Lambda)$.
\end{lemma}

{\it Proof.} $(i)$ and $(ii)$ are clear by the properties of the
conditional expectation. For $(iii)$, observe that
\begin{eqnarray*}
  \Lambda(\{{g_\Lambda}=0\})&=&\sum\limits_{e\in E}\int\limits_{\{g_\Lambda=0\}}
  1_{_1[e]}d\Lambda\\
  &=&\sum\limits_{e\in
  E}\int\limits_{\{g_{\Lambda e}=0\}}g_{\Lambda e}d\Lambda=0.
\end{eqnarray*}

For $(iv)$, by the pull-out property of the conditional expectation
( see e.g. Theorem 10.1.9 in \cite{Du}),
\begin{eqnarray*}
  \int \left(\frac{1}{{g_\Lambda}}\wedge n\right)
  d\Lambda&=&\sum\limits_{e\in E}\int
1_{_1[e]}\left(\frac{1}{g_{\Lambda e}}\wedge n\right)
  d\Lambda\\
  &=&\sum\limits_{e\in E}\int
  E_\Lambda\left(\left.1_{_1[e]}\left(\frac{1}{g_{\Lambda e}}\wedge n\right)\right|\mathcal{F}\right)
  d\Lambda\\
 &=&\sum\limits_{e\in E}\int g_{\Lambda e}\left(\frac{1}{g_{\Lambda e}}\wedge n\right)
  d\Lambda\\
  &\leq&\sum\limits_{e\in E}1,
\end{eqnarray*}
for all $n\in\mathbb{N}$. Hence, $1/{g_\Lambda}\in
\mathcal{L}^1(\Lambda)$. Since $\log x\leq x-1$,
\begin{eqnarray*}
 \int |\log g_\Lambda| d\Lambda &=&\int \log\frac{1}{g_\Lambda} d\Lambda\\
 &\leq&\int \left(\frac{1}{g_\Lambda}-1\right)d\Lambda\\
 &<&\infty.
\end{eqnarray*}
 Thus, $\log g_\Lambda\in\mathcal{L}^1(\Lambda)$.\hfill$\Box$

\begin{lemma}\label{es}
Let $\Lambda\in P_S(\Sigma)$. For each $e\in E$, let
$\dot{g}_{\Lambda e}$ be a Borel measurable version of $g_{\Lambda
e}$ such that $\log\dot{g}_{\Lambda}$ is an energy function and
$\sum_{e\in E}\dot{g}_{\Lambda e}(\sigma)\leq 1$ for all
$\sigma\in\Sigma$. Then $\Lambda$ is an equilibrium state for
 $\log\dot{g}_\Lambda$  and
\[h_\Lambda(S)+\Lambda(\log \dot{g}_\Lambda)=0,\] i.e. $P(\log
\dot{g}_\Lambda)=0$. If $\Theta\in P_S(\Sigma)$ is also an
equilibrium state for $\log \dot{g}_\Lambda$, then
$g_\Theta=\dot{g}_\Lambda$ $\Theta$-a.e..
\end{lemma}

{\it Proof.} The proof goes similarly to that of Theorem 1 in
\cite{Le}. The key for the proof is the pull-out property of the
conditional expectation ( see e.g. Theorem 10.1.9 in \cite{Du}).

Since $\mathcal{A}_0$ is a generator for the Borel $\sigma$-algebra,
we know (see e.g. \cite{W}) that
\begin{eqnarray*}
 h_\Lambda(S)=-\sum\limits_{e\in E}\int E_\Lambda\left(\left.1_{_1[e]}\right|\mathcal F\right)\log
 E_\Lambda\left(\left.1_{_1[e]}\right|\mathcal F\right)d\Lambda.
\end{eqnarray*}
By Lemma \ref{ef},  $1_{_1[e]}\log
 g_{\Lambda e}\in\mathcal{L}^1(\Lambda)$ for all $e\in E$. Therefore, by the pull-out
property of the conditional expectation,
\begin{eqnarray*}
 h_\Lambda(S)&=&-\sum\limits_{e\in E}\int 1_{_1[e]}\log
 g_{\Lambda e}d\Lambda\\
 &=&-\sum\limits_{e\in E}\int 1_{_1[e]}\log
 g_\Lambda d\Lambda\\
 &=&-\int \log g_\Lambda d\Lambda.
\end{eqnarray*}
 For the first part of the Lemma, it remains to show that
\[h_\Theta(S)+\Theta(\log\dot{g}_\Lambda)\leq 0\mbox{ for all }\Theta\in P_S(\Sigma).\]
Let $\Theta\in P_S(\Sigma)$. By the above, $h_\Theta(S)=-\Theta(\log
g_\Theta)$. If $\Theta(\{\dot{g}_\Lambda=0\})>0$, then
$h_\Theta(S)+\Theta(\log \dot{g}_\Lambda)=-\infty<0$. Otherwise,
since $\log x\leq x-1$ for all $x>0$, it follows that
\begin{eqnarray*}
 h_\Theta(S)+\Theta(\log
 \dot{g}_\Lambda)&=&\int\log\frac{\dot{g}_\Lambda}{g_\Theta}d\Theta\\
 &\leq&\int\left(\frac{\dot{g}_\Lambda}{g_\Theta}-1\right)d\Theta.
\end{eqnarray*}
By Lemma \ref{ef}, $1_{_1[e]}(\dot{g}_{\Lambda e}/g_{\Theta e}-1)\in
\mathcal{L}^1(\Theta)$ for all $e\in E$. Therefore, by the pull-out
property of the conditional expectation,
\begin{eqnarray*}
  \int\left(\frac{\dot{g}_\Lambda}{g_\Theta}-1\right)d\Theta
  &=& \sum\limits_{e\in E}\int1_{_1[e]}\left(\frac{\dot{g}_{\Lambda e}}{g_{\Theta e}}-1\right)d\Theta\\
  &=& \sum\limits_{e\in E}\int g_{\Theta e}\left(\frac{\dot{g}_{\Lambda e}}{g_{\Theta e}}-1\right)d\Theta\\
  &=& \int\sum\limits_{e\in E} \left(\dot{g}_{\Lambda e}-g_{\Theta e}\right)d\Theta\\
  &\leq& 0.
\end{eqnarray*}
Hence, \[h_\Theta(S)+\Theta(\log\dot{g}_\Lambda)\leq 0,\] i.e.
$\Lambda$ is an equilibrium state for $\log\dot{g}_\Lambda$.

Now, suppose $\Theta_0$ is also an equilibrium state for $\log\dot
{g}_\Lambda$, i.e.
\[h_{\Theta_0}(S)+\Theta_0(\log\dot{g}_\Lambda)=\sup\limits_{\Theta\in
P_S(\Sigma)}\left(h_\Theta(S)+\Theta(\log
\dot{g}_\Lambda)\right)=h_\Lambda(S)+\Lambda(\log\dot{g}_\Lambda)=0.\]

 Then, by the above, the following equality must hold true:
\[\int\log\frac{\dot{g}_\Lambda}{g_{\Theta_0}}d\Theta_0=\int\left(\frac{\dot{g}_\Lambda}{g_{\Theta_0}}-1\right)
d\Theta_0.\] But it is true if and only if
\[\log\frac{\dot{g}_\Lambda}{g_{\Theta_0}}=\left(\frac{\dot{g}_\Lambda}{g_{\Theta_0}}-1\right)\
\Theta_0\mbox{-a.e.}.\] And the latter happens if and only if
$\dot{g}_\Lambda=g_{\Theta_0}$ $\Theta_0$-a.e.. \hfill$\Box$

Now, we are going to prove what seems to be the main lemma for the
generalized Markov shift  associated with a contractive Markov
system. For that we need to define some measures on the product
space $K\times\Sigma$.

Denote by $\mathcal{A}$  the finite $\sigma$-algebra generated by
the partition $\{_0[e]:e\in E\}$ of $\Sigma$ and define, for each
integer $m\leq 1$,
\[\mathcal{A}_m:=\bigvee\limits_{i=m}^{+\infty} S^{-i}\mathcal{A},\]
which is the smallest $\sigma$-algebra containing  all finite
 $\sigma$-algebras $\bigvee_{i=m}^{n}
S^{-i}\mathcal{A}$, $n\geq m$. Let $x\in K$. For every integer
$m\leq 1$, let $P_x^m$ be a probability measure on $\sigma$-algebra
$\mathcal{A}_m$ given by
\[P^m_x( _{m}[e_{m},...,e_n])=p_{e_{m}}(x)p_{e_{m+1}}(w_{e_{m}}(x))...p_{e_n}(w_{e_{n-1}}\circ...\circ
w_{e_{m}}(x))\] for all cylinders $_{m}[e_{m},...,e_n]$, $n\geq{m}$.
By Lemma 1 from \cite{Wer3}, $x\longmapsto P_x^m(A)$ is a Borel
measurable function on $K$. Therefore, we can define, for every
integer $m\leq 0$,
\[\tilde M_m\left(A\times Q\right):=\int\limits_{A}P^m_x\left(Q\right)d\mu(x)\]
for all $A\in\mathcal{B}(K)$ and all $Q\in\mathcal{A}_m$. Then
$\tilde M_m$ extends uniquely to a probability measure on the
product $\sigma$-algebra $\mathcal{B}( K)\otimes\mathcal{A}_m$ with
\[\tilde M_m(\Omega)=\int
P^m_x\left(\left\{\sigma\in\Sigma:(x,\sigma)\in\Omega\right\}\right)d\mu(x)\]
for all $\Omega\in\mathcal{B}( K)\otimes\mathcal{A}_m$. Note that
the set of all $\Omega\in\mathcal{B}( K)\otimes\mathcal{A}_m$ for
which the integrand in the above is measurable forms a Dynkin system
which contains the set all rectangles   $A\times Q$,
$A\in\mathcal{B}(K)$, $Q\in\mathcal{A}_m$. As the latter is
$\cap$-stable and generates $\mathcal{B}(K)\otimes\mathcal{A}_m$,
the integrand is measurable for all
$\Omega\in\mathcal{B}(K)\otimes\mathcal{A}_m$. Further, note that
$P^m_x\left(\left\{\sigma\in\Sigma:(x,\sigma)\in\Omega\right\}\right)=\int
1_{\Omega}(x,\sigma)dP_x^m(\sigma)$ for all
$\Omega\in\mathcal{B}(K)\otimes\mathcal{A}_m$. Therefore
\[\int s d\tilde M_m=\int\int s(x,\sigma)dP^m_x(\sigma)d\mu(x)\] for
all $\mathcal{B}(K)\otimes\mathcal{A}_m$-simple functions $s$. Now,
let $\psi$ be a $\mathcal{B}(K)\otimes\mathcal{A}_m$-measurable and
$\tilde M_m$-integrable function on $K\times\Sigma$. Then the
 usual monotone approximation of positive and negative parts of $\psi$ by simple functions
 and the B. Levi Theorem imply that
\[\int \psi d\tilde M_m=\int\int \psi(x,\sigma)dP^m_x(\sigma)d\mu(x).\]

\begin{lemma}\label{ml}
Suppose $\mathcal{M}$ is a CMS with an invariant Borel probability
measure
 $\mu$ such that $C:=\sum_{i=1}^N\int_{K_i}d(x,x_i)\ d\mu(x)<\infty$ for some $x_i\in K$, $i=1,...,N$, and $M(Y)=1$,
where $M$ is the generalized Markov measure associated with
$\mathcal{M}$ and $\mu$.
  Then
\[E_M\left(1_{ _1[e]}|\mathcal{F}\right)=p_e\circ F\mbox{ $M$-a.e. for all }e\in E.\]
\end{lemma}

{\it Proof.}
 Fix $e\in E$.
Let $\mathcal{F}_m:=\bigvee\limits_{i=0}^m S^i(\mathcal{A})$ for all
$m\in\mathbb{Z}\setminus\mathbb{N}$. Let's use the notation
$(\sigma_m,...,\sigma_0)^*$ iff $M( _m[\sigma_m,...,\sigma_0])>0$.
Then obviously
\[E_M\left(1_{ _1[e]}|\mathcal{F}_m\right)(\tilde\sigma)=\sum\limits_{(\sigma_m,...,\sigma_0)^*}\frac{\int\limits_{
_m[\sigma_m,...,\sigma_0]}1_{ _1[e]}\ dM}{M(
_m[\sigma_m,...,\sigma_0])}1_{
_m[\sigma_m,...,\sigma_0]}(\tilde\sigma)\] for all $M$-a.e.
$\tilde\sigma\in\Sigma$. By Doob's Martingale Theorem,
\[E_M\left(1_{ _1[e]}|\mathcal{F}_m\right)\to E_M\left(1_{ _1[e]}|\mathcal{F}\right)\ M\mbox{-a.e.}.\]

Now, set
\[Z^x_m(\sigma):=w_{\sigma_0}\circ...\circ w_{\sigma_m}(x)\mbox{ and } Y_m(\sigma):=
w_{\sigma_0}\circ...\circ w_{\sigma_m}(x_{i(\sigma_m)})\] for all
$x\in K$, $\sigma\in\Sigma$ and $m\leq 0$. Then
\begin{eqnarray}\label{e}
 &&\left|E_M\left(1_{ _1[e]}|\mathcal{F}_m\right)-p_e\circ F\right|\\
&\leq&
\sum\limits_{(\sigma_m,...,\sigma_0)^*}\left|\frac{\int\limits_{K\times
_m[\sigma_m,...,\sigma_0]}p_e\circ Z^x_m(\tilde\sigma) \ d\tilde
M_m(x,\tilde\sigma)}{\tilde M_m(K\times
_m[\sigma_m,...,\sigma_0])}-p_e\circ Y_m\right|1_{
_m[\sigma_m,...\sigma_0]}+
\left|p_e\circ Y_m-p_e\circ F\right|\nonumber\\
&\leq&
\sum\limits_{(\sigma_m,...,\sigma_0)^*}\frac{\int\limits_{K\times
_m[\sigma_m,...,\sigma_0]} \left|p_e\circ
Z^x_m(\tilde\sigma)-p_e\circ Y_m\right| \ d\tilde
M_m(x,\tilde\sigma)}{\tilde M_m(K\times
_m[\sigma_m,...,\sigma_0])}1_{ _m[\sigma_m,...\sigma_0]}+
\left|p_e\circ Y_m-p_e\circ F\right|\nonumber.
\end{eqnarray}
Set \[Q_m:=\left\{(x,\sigma)\in K\times\Sigma:\
d(Z^x_m(\sigma),Y_m(\sigma))>a^{\frac{-m+1}{2}}C\right\}.\] Observe
that, by the contractiveness on average condition,
\begin{eqnarray*}
 &&\int d\left(Z^x_m(\sigma),Y_m(\sigma)\right)\ d\tilde M_m(x,\sigma)\\
&=&\int\sum\limits_{\sigma_m,...,\sigma_0}p_{\sigma_m}(x)...p_{\sigma_0}(w_{\sigma_{-1}}\circ...\circ
w_{\sigma_{m}}x)
   d\left(w_{\sigma_{0}}\circ...\circ w_{\sigma_{m}}x,w_{\sigma_{0}}\circ...\circ w_{\sigma_{m}}x_{i(\sigma_m)}\right)\ d\mu(x)\\
&\leq& a^{-m+1}\sum\limits_{i=1}^N\int\limits_{K_i}d(x,x_i)\ d\mu(x)\\
&=& a^{-m+1}C.
\end{eqnarray*}
Hence,
\[a^{\frac{-m+1}{2}}C\tilde M_m\left(Q_m\right)\leq a^{-m+1}C,\]
that is,
\[\tilde M_m\left(Q_m\right)\leq a^{\frac{-m+1}{2}}.\]
Therefore,
\begin{eqnarray}\label{i}
(\ref{e})\leq
X_m+\phi\left(a^{\frac{-m+1}{2}}C\right)+\left|p_e\circ Y_m-p_e\circ
F\right|\ M\mbox{-a.e.},
\end{eqnarray}
where \[X_m:=\sum\limits_{(e_m,...,e_0)^*}\frac{\tilde
M_m\left(Q_m\cap (K\times _m[e_m,...,e_0])\right)} {\tilde
M_m(K\times _m[e_m,...,e_0])}1_{ _m[e_m,...e_0]},\] and $\phi$ is
the modulus of uniform continuity of $p_e|_{K_{i(e)}}$. Then
\[\int X_m\ dM=\tilde M_m\left(Q_m\right)\leq a^{\frac{-m+1}{2}}.\]
Set
\[\Omega_m:=\left\{\sigma\in\Sigma:\ X_m(\sigma)>a^{\frac{-m+1}{4}}\right\}.\]
Then
\[a^{\frac{-m+1}{4}}M\left(\Omega_m\right)\leq\int\limits_{\Omega_m}X_m\ dM\leq a^{\frac{-m+1}{2}}.\]
Hence
\[M\left(\Omega_m\right)\leq a^{\frac{-m+1}{4}}.\]
Set
\[\Omega:=\bigcap\limits_{n\leq 0}\bigcup\limits_{m\leq n}\Omega_m.\]
Then $X_m(\sigma)\to 0$ for all $\sigma\in\Sigma\setminus\Omega$,
and
\[M(\Omega)=0,\] by the Borel-Cantelli argument. Since each $p_e|_{K_{i(e)}}$ is uniformly continuous and $M(Y)=1$ implies that
$Y_m\to F$ $M$-a.e., we conclude, by (\ref{i}), that
\[\left|E_M\left(1_{ _1[e]}|\mathcal{F}_m\right)-p_e\circ F\right|\to 0\; M\mbox{-a.e.}.\]
\hfill$\Box$

\begin{prop}\label{esr}
Suppose $\mathcal{M}$ is a CMS with an invariant Borel probability
measure
 $\mu$ such that $\sum_{i=1}^N\int_{K_i}d(x,x_i)\ d\mu(x)<\infty$ for some $x_i\in K$, $i=1,...,N$, and $M(Y)=1$,
where $M$ is the generalized Markov measure associated with
$\mathcal{M}$ and $\mu$. Then the following hold.\\
(i) $M$ is an equilibrium state for $u$.\\
(ii) $P(u)=0$.
\end{prop}
{\it Proof.} By Lemma \ref{ml},
\[E_M\left(1_{ _1[e]}|\mathcal{F}\right)=p_e\circ F\ M\mbox{-a.e. for all }e\in E.\]
For each $e\in E$, set
\begin{equation*}\label{ef}
    \dot{g}_{M e}(\sigma):=\left\{\begin{array}{cc}
     p_{e}\circ F(\sigma)&  \mbox{if }\sigma\in Y \\
    0& \mbox{if }\sigma\in\Sigma\setminus Y.
  \end{array}\right.
\end{equation*}
Since $M(Y)=1$, each $\dot{g}_{Me}$ is a version of $g_{Me}$ which
satisfies the hypothesis of Lemma \ref{es} and $u=\log\dot{g}_M$,
where $\dot{g}_M(\sigma):=\dot{g}_{M{\sigma_1}}(\sigma)$ for all
$\sigma\in\Sigma$. Hence, by Lemma \ref{es}, $M$ is an equilibrium
state for $u$, and also holds $(ii)$.\hfill$\Box$

\begin{prop}\label{esp}
 Let $\Lambda\in P_S(\Sigma)$ be an equilibrium state for $u$. Then
\[U^*F(\Lambda)=F(\Lambda).\]
\end{prop}
{\it Proof.} Since $\Lambda$ is an equilibrium state for $u$,
$\Lambda(Y)=1$ (otherwise $h_\Lambda(S)+\Lambda(u)=-\infty<0$).
Furthermore, by Lemma 5, $g_\Lambda=\exp u$ $\Lambda$-a.e.. Hence,
we can assume without loss of generality that
$g_\Lambda(\sigma)=\exp u(\sigma)$ for all $\sigma\in Y$.

Now, let $f\in C_B(K)$ and $_{-1}[e_0,e]\cap Y\neq\emptyset$. Then,
by  the shift invariance of $\Lambda$ and the pull-out property of
the conditional expectation,
\begin{eqnarray*}
  \int\limits_{_{-1}[e_0,e]}f\circ F\ d\Lambda&=&\int\limits_{_0[e_0]}
  1_{_1[e]}f\circ F\circ S\ d\Lambda=\int\limits_{_0[e_0]}
  1_{_1[e]}f\circ w_e\circ F\ d\Lambda\\
  &=&\int\limits_{_0[e_0]} g_{\Lambda e}f\circ w_e\circ F\
  d\Lambda.
\end{eqnarray*}
Let $\sigma\in _0[e_0]\cap Y$. Since $_{0}[e_0,e]\cap
Y\neq\emptyset$, there exists $\sigma'\in Y$ such that
$\sigma'=(...,\sigma_{-2},\sigma_{-1},e_0,e,\sigma'_2,\sigma'_3,...)$.
Hence $g_{\Lambda e}(\sigma)=g_\Lambda(\sigma')=\exp
u(\sigma')=p_e\circ F(\sigma)$. Thus
\[g_{\Lambda e}(\sigma)=p_e\circ F(\sigma)\mbox{ for all }\sigma\in\ _0[e_0]\cap Y.\]
Therefore
\begin{eqnarray*}
  \int\limits_{_0[e_0]} g_{\Lambda e}f\circ w_e\circ F\
  d\Lambda&=&\int\limits_{_0[e_0]} p_e\circ Ff\circ w_e\circ F\
  d\Lambda.
\end{eqnarray*}
Summing for all $e_0\in E$ gives
\[\int\limits_{_{0}[e]}f\circ F\ d\Lambda=\int\limits p_e\circ Ff\circ w_e\circ F\
  d\Lambda.\]
  Hence,
  \[F(\Lambda)(f)=\int f\circ F\ d\Lambda=\sum\limits_{e\in
  E}\int p_ef\circ w_e\ dF(\Lambda)=U^*F(\Lambda)(f)\]
as desired.\hfill$\Box$

\begin{theo}\label{gMme}
 Suppose CMS $\mathcal{M}$ has a unique invariant Borel probability measure $\mu$ and $M(Y)=1$, where $M$ is the associated
generalized Markov measure.
Then the following hold.\\
$(i)$ $M$ is a unique equilibrium state for the energy function $u$,\\
$(ii)$ $F(M)=\mu$,\\
$(iii)$ $h_M(S)=-\sum_{e\in E}\int_{K_{i(e)}}p_e\log p_e\ d\mu$.
\end{theo}

{\it Proof.}  By Theorem \ref{ft} (iv),
$\sum_{i=1}^N\int_{K_i}d(x,x_i)d\mu(x)<\infty$
 for all $x_i\in K_i,$ $i=1,...,N$. Therefore, by Lemma \ref{ml},
\[E_M\left(\left.1_{_1[e]}\right|\mathcal{F}\right)=p_e\circ F\ M\mbox{-a.e. for
all }e\in E.\] Hence, by Proposition \ref{esr}, $M$ is an
equilibrium state for $u$.

Now, suppose $\Lambda$ is another equilibrium state for $u$. This
implies that $\Lambda(Y)=1$ (otherwise
$h_\Lambda(S)+\Lambda(u)=-\infty<0$). Also, by Lemma \ref{es},
$g_\Lambda(\sigma)=\dot{g}_M(\sigma)$ for $\Lambda$-a.e.
$\sigma\in\Sigma$. Hence, we can assume, without loss of generality,
that
\[g_{\Lambda\sigma_1}(\sigma)=p_{\sigma_1}\circ F(\sigma)\mbox{
for all }\sigma\in Y.\] Let $[e_1,...,e_n]\subset\Sigma$ be a
cylinder set such that $[e_1,...,e_n]\cap Y\neq\emptyset$.  By the
shift-invariance of $\Lambda$,
\begin{eqnarray*}
   \Lambda([e_1,...,e_n])
   =\int\limits_{_{-n+2}[e_1,...,e_{n-1}]}1_{_1[e_n]}d\Lambda
   =
   \int\limits_{_{-n+2}[e_1,...,e_{n-1}]}g_{\Lambda{e_n}}d\Lambda.
\end{eqnarray*}
Let $\sigma\in\ _{-n+2}[e_1,...,e_{n-1}]\cap Y$. Then there exists
$\sigma'\in Y$ such that
$\sigma'=(...,\sigma_{-1},\sigma_0,e_n,\sigma'_2,\sigma'_3,...)$.
Hence, $g_{\Lambda
e_n}(\sigma)=g_{\Lambda\sigma'_1}(\sigma')=p_{\sigma'_1}\circ
F(\sigma')=p_{e_n}\circ F(\sigma)$. We conclude that
\[g_{\Lambda e_n}(\sigma)=p_{e_n}\circ F(\sigma)\mbox{
for all }\sigma\in\ _{-n+2}[e_1,...,e_{n-1}]\cap Y.\] Note that $F$
is $\mathcal{F}$-measurable and
  $F(S\sigma)=w_{\sigma_1}(F(\sigma))$ for all $\sigma\in Y$. Therefore,
\begin{eqnarray*}
   && \int\limits_{_{-n+2}[e_1,...,e_{n-1}]}g_{\Lambda{e_n}}d\Lambda\\
  &=& \int\limits_{_{-n+2}[e_1,...,e_{n-1}]}p_{e_n}\circ Fd\Lambda\\
  &=& \int\limits_{_{-n+3}[e_1,...,e_{n-2}]}1_{_1[e_{n-1}]}p_{e_n}\circ
  F\circ Sd\Lambda\\
  &=& \int\limits_{_{-n+3}[e_1,...,e_{n-2}]}1_{_1[e_{n-1}]}p_{e_n}\circ w_{e_{n-1}}\circ F
      d\Lambda.
\end{eqnarray*}
By using the pull-out property of the conditional expectation and
repeating the above argumentation, we obtain that
\begin{eqnarray*}
  && \int\limits_{_{-n+3}[e_1,...,e_{n-2}]}1_{_1[e_{n-1}]}p_{e_n}\circ w_{e_{n-1}}\circ F
      d\Lambda\\
  &=& \int\limits_{_{-n+3}[e_1,...,e_{n-2}]}p_{e_{n-1}}\circ Fp_{e_n}\circ w_{e_{n-1}}\circ F
      d\Lambda\\
  &.&\\
  &.&\\
  &.&\\
  &=& \int p_{e_1}\circ F p_{e_2}\circ w_{e_1}\circ F...p_{e_n}\circ w_{e_{n-1}}\circ...\circ w_{e_1}\circ
  Fd\Lambda\\
  &=& \int p_{e_1} p_{e_2}\circ w_{e_1}...p_{e_n}\circ w_{e_{n-1}}\circ...\circ
  w_{e_1}dF(\Lambda).
\end{eqnarray*}

 Thus,
the equality $\Lambda=M$ will follow from $F(\Lambda)=\mu$, but this
follows
 by the uniqueness of the invariant measure $\mu$, since
 $U^*F(\Lambda)=F(\Lambda)$ by Proposition \ref{esp}.
Thus, the claims $(i)$ and $(ii)$ hold true. By Proposition
\ref{esr} $(ii)$ and Lemma \ref{ml},
\begin{eqnarray*}
h_M(S)&=&-\sum\limits_{e\in E}\int 1_{ _1[e]}\log p_e\circ F\ dM\\
        &=&-\sum\limits_{e\in E}\int p_e\circ F\log p_e\circ F\ dM\\
        &=&-\sum\limits_{e\in E}\int\limits_{K_{i(e)}} p_e\log p_e\ d\mu.
\end{eqnarray*} This proves $(iii)$.
\hfill$\Box$

\begin{cor}\label{DCc}
Suppose $\left(K_{i(e)},w_e,p_e\right)_{e\in E}$ is an irreducible
CMS such that each $p_e|_{K_{i(e)}}$ is Dini-continuous and bounded
away from zero.
Then the following hold.\\
$(i)$ The generalized Markov measure $M$ is a unique equilibrium state for the energy function $u$,\\
$(ii)$ $P(u)=0$,\\
$(iii)$ $F(M)=\mu$,\\
$(iv)$ $h_M(S)=-\sum_{e\in E}\int_{K_{i(e)}}p_e\log p_e\ d\mu$.
\end{cor}
{\it Proof.} By Theorem 2 in \cite{Wer1}, the CMS has a unique
invariant Borel probability measure. Since $M(Y)=1$ by Lemma
\ref{pY} $(ii)$, the claims follow by Theorem
\ref{gMme}.\hfill$\Box$

\begin{Remark}
 The author would like to point out that a similar entropy formula as that proved in Theorem \ref{gMme} (iv)
 plays a central role in the recent book of Wojciech Slomczynski \cite{S}.
\end{Remark}

Finally, we would like to make some remarks on why the result
presented here might be interesting for the general theory of
thermodynamic formalism.

\begin{Remark}\label{ctf}
 First of all, recall that the theory of equilibrium states
 is presented usually only for upper semicontinuous energy
 functions (see e.g. \cite{Kel}). The uniqueness of an equilibrium state is
 known  on  sub-shifts of finite type in general only for energy functions satisfying some stronger continuity
 conditions, e.g. Dini-continuity (note that the Dini-continuity and the {\it regularity} of a
function coincide on  one-dimensional lattices) \cite{Kel},
\cite{Ru}.

If all maps $w_e$ of $\mathcal{M}$ are contractive and all
probabilities $p_e|_{K_{i(e)}}$ are Dini-continuous, then the coding
map $F$ is defined everywhere on $\Sigma_G$ and is
H\"{o}lder-continuous. Hence, the energy function $u$ is
upper-semicontinuous and $u|_{\Sigma_G}$ is Dini-continuous (easy to
check, since
 $\log x\leq x-1$). In this case, Corollary
\ref{DCc} (i) fits nicely into the well known thermodynamic
formalism.

Now, let us consider $u$ if $w_e$'s are contractive only on average.
In this case, $Y$ is not necessarily
 closed (see Example \ref{ex}), i.e. $u$ is not necessarily upper semicontinuous. Therefore, even the
existence of an equilibrium state for $u$ is not guaranteed by the
existing thermodynamic formalism. Moreover, by Lemma \ref{cm}, we
only know that for all $\epsilon>0$ there exists $Q\subset Y$ with
$M(Q)>1-\epsilon$ such that $u|_Q$ is Dini-continuous, but the sum
$\sum_{k=1}^\infty\phi_Q(2^{-k})$, where $\phi_Q$ is the modulus of
uniform continuity of $u|_Q$, increases if we choose $Q$ larger.

Summing up, Corollary \ref{DCc} shows that the general contractive
Markov systems considered here still inherit some of their
thermodynamic properties from finite Markov chains, even though
their energy function belongs to a class which, as far as the author
is aware, is not considered by the existing theory of thermodynamic
formalism.
\end{Remark}

\begin{Remark}
An important result of the thermodynamic formalism on topologically
mixing
 subshifts of finite type is that for a
Dini-continuous energy function the unique equilibrium state can be
obtained as a unique {\it Gibbs state} for the same energy function
(see e.g. \cite{Kel}, \cite{Ru}).

If all $w_e$'s are contractive, the energy function $u|_{\Sigma_G}$
is Dini-continuous  (as in Remark \ref{ctf}) and therefore, by
Corollary \ref{DCc} (i), $M$ is also a unique Gibbs state for $u$.
However, $Y$ is not necessarily a subshift of finite type (not
necessarily closed) and $u|_Y$ is not necessarily Dini-continuous if
the maps are contractive only on average. We do not know in this
case whether the measure $M$ still can be constructed as a unique
Gibbs state.
\end{Remark}

\subsection*{Acknowledgements}
I would like to thank: EPSRC and School of Mathematics
 and Statistics of University of St Andrews for providing me with a scholarship and excellent working
 conditions in St Andrews, the anonymous referees for suggestions on improvements for this paper. Also,
I would like to thank Barry Ridge and Wang Yang for their help in
the production of this paper.

\end{document}